\documentclass{birkmult}
\usepackage{epsfig}
\usepackage{amssymb}
\usepackage{amsmath}
\usepackage{amsfonts}

\newcommand{\kt}{\rangle}
\newcommand{\bra}{\langle}
\newcommand{\Sc}{{\mathcal S}}
\newcommand{\A}{{\mathfrak A}}
\newcommand{\Hil}{\mathcal H}
\newcommand{\1}{1 \!\! 1}
\newcommand{\Lc}{{\mathcal L}}

\def\be{\begin{equation}}
\def\ee{\end{equation}}
\def\ba{\begin{array}{c}}
\def\ea{\end{array}}

 \theoremstyle{definition}
 
 \theoremstyle{remark}

 \numberwithin{equation}{section}

\begin{document}
%
%
%
%
%
%
%
%
%
\title[]{The dynamical problem for a non self-adjoint Hamiltonian}

\author{Fabio Bagarello}

\address{%
Dipartimento di Metodi e Modelli Matematici, Facolt\`a di Ingegneria,
Universit\`a di Palermo, I-90128
Palermo, Italy.
}
\email{bagarell@unipa.it}

\author{Miloslav Znojil}

\address{%
 Nuclear Physics Institute ASCR\br
 250 68 \v{R}e\v{z}, Czech Republic
}
\email{znojil@ujf.cas.cz}
\begin{abstract}

We review some recent results of the so-called quasi-hermitian
quantum mechanics, with particular focus on the quantum dynamics
both in the Schr\"odinger and in the Heisenberg representations. The
role of Krein spaces is also discussed.

\end{abstract}
\maketitle

\section{Introduction}

In the analysis of the dynamics of a closed quantum system ${\Sc}$ a
special role is played by the energy  $H$, which is typically the
self-adjoint operator defined by the sum of the kinetic energy of
$\Sc$ and of the potential energy giving rise to the conservative
forces acting on $\Sc$. The {\em most common} approaches in the
description of $\Sc$ are the following:

\begin{enumerate}

\item
{\em the algebraic description (AD)}: in this approach the
observables of $\Sc$ are elements of a C*-algebra $\A$
(which coincides with  $B(\Hil)$ for some Hilbert space $\Hil$).
This means, first of all, that  $\A$ is a vector space over
$\mathbb{C}$ with a multiplication law such that $\forall A,B\in\A$,
$AB\in\A$. Also, two such elements can be summed up and the
following properties hold: $\forall A,B,C\in\A$ and $\forall
\alpha,\beta\in\mathbb{C}$ we have
$$
A(BC)=(AB)C,\hspace{5mm} A(B+C)=AB+AC, \hspace{5mm} (\alpha A)(\beta
B)=\alpha\beta (AB).
$$
An involution is a map $*:\A\rightarrow \A$ such that
$$
A^{**}=A,\hspace{5mm} (AB)^*=B^*A^*,\hspace{5mm} (\alpha A+\beta
B)^*=\overline{\alpha}\,A^*+\overline{\beta}\,B^*
$$
A *-algebra $\A$ is an algebra with an involution *. $\A$ is  a
{\it normed algebra} if there exists a map, {\it the norm of the
algebra}, $\|.\|:\A\rightarrow \mathbb{R}_+$, such that:
$$
\|A\|\geq 0,\hspace{5mm} \|A\|=0 \Longleftrightarrow
A=0,\hspace{5mm} \|\alpha A\|=|\alpha|\,\|A\|,$$ $$\|A+B\|\leq
\|A\|+\|B\|,\hspace{5mm}\|AB\|\leq \|A\|\,\|B\|.
$$

If $\A$ is complete wrt $\|.\|$, then it is called a {\it   Banach
algebra}, or a {\it   Banach *-algebra} if $\|A^*\|=\|A\|$. If
further $\|A^*A\|=\|A\|^2$ holds for all $A\in\A$, then $\A$ is a
{\it   C*-algebra}.

The {\it states} are linear, positive and normalized functionals on
$\A$, which look like $\rho(\hat A)=tr(\hat\rho A)$, where
$\A=B(\Hil)$,   $\hat\rho$ is a trace-class operator and $tr$ is the
trace on $\Hil$. This means in particular that
$$
\rho(\alpha_1 A+\alpha_2 B)=\alpha_1\rho(A)+\alpha_2\rho(B)$$ and
that, if $\A$ has the identity $\1$,
$$
\rho(A^*A)\geq 0;\hspace{5mm} \rho(\1)=1.
$$
An immediate consequence of these assumptions, and in particular of
the positivity of $\rho$, is that $\rho$ is also continuous, i.e.
that $|\rho(A)|\leq \|A\|$ for all $A\in \A$.

The dynamics in the Heisenberg representation for the closed
quantum system $\Sc$ is given by the map

$$\A\ni A\rightarrow \alpha^t(A)=U_tAU_t^\dagger\in \A,\,\forall
t$$ which defines a 1-parameter group of *-automorphisms of $\A$
satisfying the following conditions
$$\alpha^t(\lambda A)=\lambda \alpha^t(A),\hspace{5mm}
\alpha^t(A+B)=\alpha^t(A)+\alpha^t(B),$$
$$\alpha^t(AB)=\alpha^t(A)\,\alpha^t(B),\hspace{5mm}
\|\alpha^t(A)\|=\|A\|,\, \mbox{ and }\,
\alpha^{t+s}=\alpha^t\,\alpha^s.$$

In the Schr\"odinger representation the time evolution is the dual
of the one above, i.e. it is the  map between states defined by $\hat
\rho\rightarrow \hat\rho_t={\alpha^t}^*\hat\rho$.

\item
{\em the Hilbert space description (HSD)}: this is much simpler, at
a first sight. We work in some fixed Hilbert space $\Hil$, somehow
related to the system we are willing to describe, and we proceed as
follows:

each observable $A$ of the physical system corresponds to a
self-adjoint operator $\hat A$ in $\Hil$;

the pure states of the physical system corresponds to normalized
vectors of $\Hil$;

the expectation values of $A$ correspond to the following mean
values: $<\psi,\hat A\psi>=\rho_\psi(\hat A)=tr(P_\psi\,\hat A)$,
where we have also introduced a projector operator $P_\psi$ on
$\psi$ and $tr$ is the trace on $\Hil$;

the states which are not pure, i.e. the mixed states, correspond to
convex linear combinations $\hat\rho=\sum_j\,w_j\,\rho_{\psi_n}$,
with $\sum_j\,w_j=1$ and $w_j\geq 0$ for all $j$;

the dynamics (in the Schr\"odinger representation) is given by a
unitary operator $U_t:=e^{i H t/\hbar}$, where $H$ is the
self-adjoint energy operator, as follows: $\hat\rho \rightarrow
\hat\rho_t=U_t^\dagger\hat\rho U_t$. In the Heisenberg
representation the states do not evolve in time while the operators
do, following the  {\it dual} rule:
 $\hat A \rightarrow \hat A_t=U_t\hat A U_t^\dagger$, and the Heisenberg
 equation of motion is satisfied:
$\frac{d}{dt}\hat A_t=\frac{i}{\hbar}[H,\hat A_t]$. It is very well
known that these two different representations have the same
physical content: indeed we have $\hat\rho(\hat A_t)=\hat\rho_t(\hat
A)$, which means that what we measure in experiments, that is the
time evolution of the mean values of the observables of $\Sc$, do
not depend on the representation chosen.

\end{enumerate}

The $AD$ is especially useful when $\Sc$ has an infinite number of
degrees of freedom, \cite{bra, sew1, bagrev}, while the $HSD$ is
quite common for ordinary quantum mechanical systems, i.e. for those
systems with a finite number of degrees of freedom. The reason why
the algebraic approach to ordinary quantum mechanics is not very
much used in this simpler case follows from the following von
Neumann uniqueness theorem: {\it for finite quantum mechanical
systems there exists only one irreducible representation}. This
result is false for systems with infinite degrees of freedom
(briefly, in $QM_\infty$), for which $AD$ proved to be  useful, for
instance in the description of phase-transitions \cite{thirr}.

As we have already said, in the most common applications of quantum
theory the self-adjoint Hamiltonian is just the sum of a kinetic
plus an interaction term, and the Hilbert space in which the model
is described is usually $\Hil:=\Lc^2(\mathbb{R}^D)$, for $D=1, 2$ or
3. Of course any unitary map defined from $\Hil$ to any (in general
different) other Hilbert space $\tilde\Hil$ does not change the
physics which is contained in the model, but only provides a
different way to extract the results from the model itself. In
particular, the mean values of the different observables do not
depend from the Hilbert space chosen, as far as the different
representations are connected by unitary maps. These observables are
other self-adjoint operators whose eigenvalues are interesting for
us since they have some physical meaning. For a quantum particle
moving in $D-$dimensional Euclidian space, for example, people
usually work with the position operator $\mathfrak{q} $, whose
eigenvalues are used to label the wave function $\psi(\vec{x},t)$ of
the system, which is clearly an element of $\Lc^2(\mathbb{R}^D)$. It
might be worth reminding that we are talking here of {\em
representations} from two different points of view: the Heisenberg
and the Schr\"odinger representations are two (physical) equivalent
ways to describe the dynamics of $\Sc$, while in the $AD$ a
representation is a map from $\A$ to $B(\Hil)$ for a chosen $\Hil$
which preserve the algebraic structure of $\A$. Not all these kind
of representations are unitarily equivalent, and for this reason
they can describe different physics (e.g. different phases of
$\Sc$), \cite{bra,sew1,bagrev}.

Recently, Bender and Boettcher \cite{BB} emphasized that many
Hamiltonians $H$ which look unphysical in  $\Hil$ may still be
correct and physical, provided only that the conservative textbook
paradigm is replaced by a modification called ${\bf PT}-$symmetric
quantum mechanics (PTSQM, cf., e.g., reviews
\cite{Dorey,Carl,ali,SIGMA} for more details). Within the PTSQM
formalism the operators ${\bf P}$ and  ${\bf T}$ (which characterize
a symmetry of the quantum system in question) are usually
pre-selected as parity and time reversal, respectively.

\section{Quantization recipes using non-unitary Dyson mappings $\Omega$}




The main appeal of the PTSQM formalism lies in the permission of
selecting, for  phenomenological purposes, various new and
nonstandard Hamiltonians exhibiting the manifest non-hermiticity
property $H \neq  H^\dagger$ in $\Hil$. It is clear that these
operators cannot generate unitary time evolution in $\Hil$ via
exponentiation, \cite{bagexp}, but this does not exclude
\cite{Geyer} the possibility of finding an unitary time evolution in
a different Hilbert space, not necessarily uniquely determined
\cite{sqw}, which, following the notation in \cite{SIGMA}, we
indicate as $\Hil^{(P)}$, $P$ standing for {\em physical}. Of
course, the descriptions of the dynamics in $\Hil$ and $\Hil^{(P)}$
cannot be connected by a unitary operator, but still other
possibilities are allowed and, indeed, these different choices are
those relevant for us here.

\vspace{2mm}

In one of the oldest applications of certain specific Hamiltonians
with the property $\hat{H} \neq \hat{H}^\dagger$ in $\Hil$ in the so
called interacting boson models of nuclei \cite{Geyer} it has been
emphasized that besides the above-mentioned fact that each such
Hamiltonian  admits many non-equivalent physical interpretations
realized via mutually nonequivalent Hilbert spaces $\Hil^{(P)}$, one
can also start from a {\em fixed} self-adjoint Hamiltonian
$\mathfrak{h}$ defined in $\Hil^{(P)}$ and move towards {\em many}
alternative isospectral images defined,  in some (different) Hilbert
space $\Hil$, by formula $\hat{H} =\Omega^{-1}\,\mathfrak{h}
\,\Omega$. Here the so called Dyson map $\Omega$ should be assumed
nontrivial i.e., non-unitary: $\Omega^\dagger\Omega=\Theta\neq \1$.

In phenomenology and practice, the only reason for preference and
choice between $\hat{H}$ and $\mathfrak{h}$ is the feasibility of
calculations and the constructive nature of  experimental
predictions. However, we also should be aware of the fact that
$\Omega$ is rather often an unbounded operator, so that many
mathematical subtle points usually arise when moving from
$\mathfrak{h}$ to $\hat H$ in the way suggested above because of,
among others, {\em domain details}. Examples of this kind of problems
are discussed in \cite{bagpb1,bagpb2,bagpb3} in connection with the
so-called {\em pseudo-bosons}.

\subsection{The three-space scenario}

For certain  complicated quantized systems (say, of tens or hundreds
of fermions as occur, typically, in nuclear physics, or for
many-body systems) the traditional theory forces us to work with an
almost prohibitively complicated Hilbert space $\Hil^{(P)}$ which is
not ``friendly" at all. Typically, this space acquires the form of a
multiple product $\bigotimes L^2(\mathbb{R}^3)$ or, even worse, of
an antisymmetrized Fock space. In such a case, unfortunately, wave
functions $\psi^{(P)}(t)$ in $\Hil^{(P)}$ become hardly accessible
to explicit construction. The technical difficulties make
Schr\"{o}dinger's equation practically useless: no time evolution
can be easily deduced.

A sophisticated way towards a constructive analysis of similar
quantum systems has been described by Scholtz et al \cite{Geyer}.
They felt inspired by the encouraging practical experience with the
so called Dyson's mappings $\Omega^{(Dyson)}$ between bosons and
fermions in nuclear physics. Still, their technique of an efficient
simplification of the theory is independent of any particular
implementation details. One only has to assume that the {\em
overcomplicated} realistic Hilbert space $\Hil^{(P)}$ is being
mapped on a {\em much simpler}, friendly intermediate space
$\Hil^{(F)}$. The latter space remains just auxiliary and unphysical
but it renders the calculations (e.g., of spectra) feasible. In
particular, the complicated state vectors $\psi^{(P)}(t)$ are made
friendlier via an invertible transition from $\Hil^{(P)}$ to
$\Hil^{(F)}$,
 \be
 \psi^{(P)}(t)=\Omega\,\psi(t)\ \ \in \ \ \Hil^{(P)}\,,
 \ \ \ \ \psi(t)\ \ \in \ \ \Hil=\Hil^{(F)}\,.
 \label{mapping}
 \ee
The introduction of the redundant superscript $^{(F)}$ underlines
the maximal friendliness of the space (note, e.g., that in the
above-mentioned nuclear-physics context of Ref.~\cite{Geyer}, the
auxiliary Hilbert space $\Hil^{(F)}$ was a bosonic space). It is
clear that, since $\Omega$ is not unitary, the inner product between
two functions $\psi_1^{(P)}(t)$ and $ \psi_2^{(P)}(t)$ (treated as
elements of  $\Hil^{(P)}$) differs from the one between $\psi_1(t)$
and $ \psi_2(t)$,
 \be
 \bra \psi_1^{(P)},\psi_2^{(P)}\kt_P=
  \bra \psi_1,\Omega^\dagger\,\Omega\,\psi_2\kt_F\ \neq \
  \bra \psi_1,\psi_2\kt_F\,.
  \label{nerovrov}
 \ee
Here we use the suffixes $P$ and $F$ to stress the fact that the
Hilbert spaces $\Hil^{(F)}$ and $\Hil^{(P)}$ are different and,
consequently, they have different inner products in general. Under
the assumption that $\ker\{\Omega\}$ only contains the zero vector,
and assuming for the moment that $\Omega$ is bounded,
$\Theta:=\Omega^\dagger\,\Omega$ may be interpreted as producing
another, alternative inner product between elements $\psi_1(t)$ and
$ \psi_2(t)$ of $\Hil^{(F)}$. This is because $\Theta$ is strictly
positive. This suggests to introduce a third Hilbert space,
$\Hil^{(S)}$, which coincides with $\Hil^{(F)}$ but for the inner
product. The new product, $\bra .,.\kt^{(S)}$, is  defined by
formula
 \be
 \bra \psi_1,\psi_2\kt^{(S)}:=
  \bra \psi_1,\Theta\psi_2\kt_F
  \label{rovnerov}
 \ee
and, because of (\ref{nerovrov}), exhibits the following
unitary-equivalence property
 \be
 \bra \psi_1,\psi_2\kt^{(S)}\ \equiv \
  \bra \psi_1^{(P)},\psi_2^{(P)}\kt_P\,.
  \label{rovnexxrov}
 \ee

More details on this point can be found in \cite{SIGMA}. It is worth
stressing that the fact that $\Omega$ is bounded makes it possible
to have $\Theta$ everywhere defined in $\Hil^{(F)}$. Under the more
general (and, we should say, more common) situation when $\Omega$ is
not bounded,  we should be careful about the possibility of
introducing $\Theta$, since $\Omega^\dagger\Omega$ could not be well
defined \cite{Ralph}: indeed, for some $f$ in the domain of
$\Omega$, $f\in D(\Omega)$, we could have that $\Omega f\notin
D(\Omega^\dagger)$. If this is the case, the best we can have is
that the inner product $\bra .,.\kt^{(S)}$ is defined on a dense
subspace of $\Hil^{(F)}$.

\subsection{A redefinition of the conjugation}
%

The adjoint $X^\dagger$ of a given (bounded) operator $X$ acting on
a certain Hilbert space ${\mathcal K}$, with inner product $(.,.)$,
is defined by the following equality:
$$
( X\varphi,\Psi)=(\varphi,X^\dagger\Psi).
$$
Here $\varphi$ and $\Psi$ are arbitrary vectors in ${\mathcal K}$.
It is clear that changing the inner product also produces a
different adjoint. Hence the adjoint in $\Hil^{(S)}$ is different
from that in $\Hil^{(F)}$, since their inner products are different.
The technical simplifying assumption that $X$ is bounded is ensured,
for instance, if we consider finite dimensional Hilbert spaces. In
this way we can avoid difficulties which could arise, e.g., due to
the unboundedness of the metric operator $\Theta$. The choice of
${\rm dim}\, \Hil^{(P,F,S)}<\infty$ gives also the chance of
getting analytical results which, otherwise, would be out of our
reach.

Once we have, in $\Hil^{(P)}$, the physical, i.e., safely hermitian
and self-adjoint $\mathfrak{h} =\Omega\,\hat{H}
\,\Omega^{-1}=\mathfrak{h}^\dagger$, we may easily deduce that
 \be
 \hat{H}=\Theta^{-1}
 \hat{H}^\dagger\,\Theta:=
 \hat{H}^\ddagger
 \label{cryptoh}
 \ee
where $^\dagger$ stands for the conjugation in either $\Hil^{(P)}$
or $\Hil^{(F)}$ while $^\ddagger$ may be treated as meaning the
conjugation in $\Hil^{(S)}$ which is metric-mediated (sometimes also
called ``non-Dirac conjugation" in physics literature). Any operator
$\hat H$ which satisfies Eq.~(\ref{cryptoh}) is said to be {\em
quasi-Hermitian}. The similarity of superscripts $^\dagger$ and
$^\ddagger$ emphasizes the formal parallels between the three
Hilbert spaces $\Hil^{(P,F,S)}$.

Once we temporarily return to the point of view of physics we must
emphasize that the use of the nontrivial metric $\Theta$ is strongly
motivated by the contrast between the simplicity of $\hat{H}$
(acting in friendly $\Hil^{(F)}$ as well as in the non-equivalent
but physical $\Hil^{(S)}$) and the practical intractability of its
isospectral partner $\mathfrak{h}$ (defined as acting in a
constructively inaccessible Dyson-image space $\Hil^{(P)}$).
Naturally, once we make the selection of $\Hil^{(S)}$ (in the role
of the space in which the quantum system in question is
represented), all of the other operators of observables (say,
$\hat{\Lambda}$) acting in $\Hil^{(S)}$ must be also self-adjoint in
the same space, i.e., they must be quasi-Hermitian with respect to
{\em the same} metric,
 \be
 \hat{\Lambda}=
 \hat{\Lambda}^\ddagger:=\Theta^{-1}
 \hat{\Lambda}^\dagger\,\Theta\,.
 \label{cryptolam}
 \ee
In opposite direction, once we start from a given Hamiltonian
$\hat{H}$ and search for a metric $\Theta$ which would make it
quasi-Hermitian (i.e., compatible with the requirement
(\ref{cryptoh})), we reveal that there exist {\em many different}
eligible metrics $\Theta=\Theta(\hat{H})$. In such a situation the
simultaneous requirement of the quasi-Hermiticity of another
operator imposes new constraints of the form
$\Theta(\hat{H})=\Theta(\hat{\Lambda})$ which restricts, in
principle, the ambiguity of the metric \cite{Geyer}. Thus, a finite
series of the quasi-Hermiticity constraints
 \be
 \hat{\Lambda}_j=
 \hat{\Lambda}^\ddagger_j\,,
 \ \ \ \ \ \ \ j = 1, 2, \ldots, J\,
 \label{cryptolamy}
 \ee
often leads to a unique physical metric $\Theta$ and to the unique,
optimal Hilbert-space representation $\Hil^{(S)}$. This is what
naturally extends the usual requirement of the ordinary textbook
quantum mechanics which requires the set of the observables to be
self-adjoint in a {\em pre-selected } Hilbert-space representation
or very concrete realization $\Hil^{(F)}$.

Let us now briefly describe some consequences of (\ref{cryptoh}) to
the dynamical analysis of $\Sc$. As we have already seen relation
(\ref{rovnerov}) defines the {\em physical} inner product. We will
now verify that this is the natural inner product to be used to find
the expected unitary evolution generated by the quasi-Hermitian
operator $\hat H$.

The first consequence of property $\hat{H}=\Theta^{-1}
\hat{H}^\dagger\,\Theta$ is that $e^{i\hat H^\dagger t}=\Theta
e^{i\hat Ht}\Theta^{-1}$. The proof of this equality is trivial
whenever the operators involved are bounded, condition which we will
assume here for simplicity (as stated above, we could simply imagine
that our Hilbert spaces are finite-dimensional). Condition
(\ref{cryptoh}) allows us to keep most of the standard approach to
the dynamics of the quantum system sketched in the Introduction,
even in presence of a manifest non-Hermiticity of $\hat H$ in
friendly $\Hil^{(F)}$. Indeed, once we consider the Schr\"odinger
evolution of a vector $\Psi(0)$, which we take to be
$\Psi(t)=e^{-i\hat Ht}\Psi(0)$, we may immediately turn attention to
the time-dependence of the {\em physical} norm in  $\Hil^{(S)}$,
$$
\|\Psi(t)\|^2:=\bra \Psi(t),\Psi(t)\kt^{(S)} =\bra e^{-i\hat
Ht}\Psi(0),\Theta e^{-i\hat Ht}\Psi(0)\kt_F=$$
$$=\bra \Psi(0),e^{i\hat H^\dagger t}\Theta e^{-i\hat Ht}
\Psi(0)\kt_F=\bra \Psi(0),\Theta \Psi(0)\kt_F=\|\Psi(0)\|^2,
$$
Thus, for all $\Psi\in\Hil^{(S)}$ the natural use of the inner
product $\bra .,.\kt^{(S)}$ and of the related norm gives a
probability which is preserved in time. This observation has also an
interesting (even if expected) consequence: the dual evolution, i.e.
the time evolution of the observables in the Heisenberg
representation, has the standard form, $X(t)=e^{i\hat Ht}Xe^{-i\hat
Ht}$, so that $\dot X(t)=i[\hat H,X(t)]$. This is true independently
of the fact that $\hat H$ is self-adjoint or not, as in the present
case. Indeed if we ask the mean value (in the product $(S)$) of the
observables to be independent of the representation adopted, i.e. if
we require that
$$
\bra \varphi(t),X\Psi(t)\kt^{(S)}
=\bra e^{-i\hat Ht}\varphi(0),\Theta Xe^{-i\hat Ht}\Psi(0)\kt
=\bra \varphi(0),X(t)\Psi(0)\kt^{(S)}
$$
then we are forced to put $X(t)=e^{i\hat Ht}Xe^{-i\hat Ht}$ (rather
than the maybe more natural $e^{i\hat Ht}Xe^{-i\hat H^\dagger t}$).
This means that a consistent approach to the dynamical problem can
be settled up also when the energy operator of the system is not
self-adjoint, paying the price by replacing the Hilbert space in
which the model was first defined, and its inner product, with
something slightly different. In this different Hilbert space the
time evolution {\em does its job}, preserving probabilities and
satisfying the usual differential equations, both in the
Schr\"odinger and in the Heisenberg pictures.

\subsection{${\bf PT}-$symmetric quantum mechanics }

One of the most efficient suppresions of the ambiguity of the metric
has been proposed within the PTSQM formalism where the role of the
additional observable $\hat{\Lambda}$ is being assigned to another
involution ${\bf C}$  \cite{Carl}. The presence as well as a
``hidden" mathematical meaning of the original involution ${\bf P}$
may be further clarified by introducing, together with our three
Hilbert spaces   $\Hil^{(P,F,S)}$, also another auxiliary space
${\bf K}$ with the structure of the Krein space endowed with an
invertible indefinite metric equal to the parity operator as
mentioned above, ${\bf P}={\bf P}^\dagger$ \cite{Azizov}. One
requires that the given Hamiltonian proves ${\bf P}-$self-adjoint in
${\bf K}$, i.e., that it satisfies equation
 \be
 \hat{H}^\dagger\,{\bf P}={\bf P}\,\hat{H}\,.
 \label{PT}
 \ee
This is an intertwining relation between two non self-adjoint
operators $\hat H$ and $\hat H^\dagger$, and $P$ is the {\em
intertwining operator} (IO, \cite{Dieudonne}). In the standard
literature on IO, see \cite{intop} for instance, the operators
intertwined by the IO are self-adjoint, so that their eigenvalues
are real, coincident, and the associated eigenvectors are orthogonal
(if the degeneracy of each eigenvalue is 1). Here, see below, these
eigenvalues are not necessarily real. Nevertheless, many situations
do exist in the literature in which the eigenvalues of some non
self-adjoint operator can be computed and turn out to be real,
\cite{bagpb1,bagpb2,bagpb3,realsp}, even if $\hat H$ is not
self-adjoint in $\Hil^{(F)}$.

Let us go back to the  requirement (\ref{PT}). Multiple examples of
its usefulness may be found scattered in the literature
\cite{realsp}. Buslaev and Grecchi \cite{BG} were probably the first
mathematical physicists who started calling property (\ref{PT}) of
the Hamiltonian a ``${\bf PT}-$symmetry". Let us now briefly explain
its mathematical benefits under a not too essential additional
assumption that the spectrum $\{E_n\}$ of $\hat{H}$ is discrete and
non-degenerate (though, naturally, not necessarily real). Then we
may solve the usual Schr\"{o}dinger equation for the (right)
eigenvectors of $\hat{H}$,
 \be
 \hat{H}\,\psi_n=E_n\,\psi_n\,,\ \ \ \ n=0,1,\ldots
 \label{SE}
 \ee
as well as its Schr\"{o}dinger-like conjugate
rearrangement
 \be
 \hat H^\dagger\,\psi^n=E_n^*\,\psi^n\,,\ \ \ \ n = 0, 1, \ldots\,.
 \ee
In the light of Krein-space rule (\ref{PT}) the latter equation
acquires the equivalent form
 \be
 \hat H\,\left ({\bf P}^{-1}
 \,\psi^n\right )=E_n^*\,\left ({\bf P}^{-1}\,\psi^n\right )
 \,,\ \ \ \ n = 0, 1,
 \ldots\,
 \label{leftSE}
 \ee
so that we may conclude that

\begin{itemize}

\item
either $E_n=E_n^*$ is real and  the action of ${\bf P}^{-1}$ on
$\psi^n$ gives merely a vector proportional to the
 $\psi_n$,

\item
or $E_n\neq E_n^*$ is not real. In this regime we have $E_n^*= E_m$
at some $m \neq n$. This means that the action of ${\bf P}^{-1}$ on
the left eigenket $\psi^n$ produces a right eigenvector of $\hat{H}$
which is proportional to certain right eigenket $\psi_m$ of
Eq.~(\ref{SE}) at a {\em different} energy, $m \neq n$.

\end{itemize}

 \noindent
We arrived at a dichotomy: Once we take all of the $n-$superscripted
left eigenvectors $\psi^n$ of $\hat{H}$ and premultiply them by the
inverse pseudometric, we obtain a new set of ket vectors
$\phi_n={\bf P}^{-1}\,\psi^n$ which are either all proportional to
their respective $n-$subscripted partners $\psi_n$ (while the whole
spectrum is real) or not. This is a way to distinguish between two
classes of Hamiltonians. Within the framework of quantum mechanics
only those with the former property of having real eigenvalues may
be considered physical (see \cite{bagpb2} for a typical
illustration). In parallel, examples with the latter property may be
still found interesting beyond the limits of quantum mechanics,
i.e., typically, in classical optics \cite{Makris}. The recent
growth of interest in the latter models (exhibiting the so called
spontaneously broken ${\bf PT}-$symmetry) may be well documented by
their presentation via the dedicated webpages~\cite{Hook}.

The Hamiltonians  $\hat{H}$ with the real and non-degenerate spectra
may be declared acceptable in quantum mechanics. In the subsequent
step our knowledge of the eigevectors $\psi^n$ of $\hat{H}$ enables
us to define the positive-definite operator of the metric directly
\cite{Ali},
 \be
 \Theta =\sum_{n=0}^\infty\,|\psi^n\kt\,\bra \psi^n|\,.
 \label{spect}
 \ee
As long as this formula defines different matrices for different
normalizations of vectors $|\psi^n\kt$ (cf. \cite{SIGMAold} for
details) we may finally eliminate this ambiguity via the
factorization ansatz
 \be
 \Theta={\bf PC}
 \label{fakto}
 \ee
followed by the {\em double} involutivity constraint \cite{Carl}
 \be
 {\bf P}^2=I\,,\ \ \ \ \ {\bf C}^2=I\,.
 \label{faktog}
 \ee
We should mention that the series in (\ref{spect}) could be just
formal. This happens whenever the operator $\Theta$ is not bounded.
This aspect was considered in many details in connection with
pseudo-bosons, see \cite{bagpb1} for instance, where one of us has
proved that $\Theta$ being bounded is equivalent to $\{\Psi^n\}$
being a Riesz basis.

Skipping further technical details we may now summarize: The operators
of the form (\ref{fakto}) +  (\ref{faktog}) have to satisfy a number
of additional mathematical conditions before they may be declared
the admissible metric operators determining the inner product in
$\Hil^{(S)}$. {\it Vice versa}, once we satisfy these conditions
(cf., e.g., Refs.~\cite{Geyer,Petr} for their list) we may replace
the unphysical and auxiliary Hilbert space $\Hil^{(F)}$ and the
intermediate Krein space ${\bf K}$ (with its indefinite metric ${\bf
P}={\bf P}^\dagger$)  by the ultimate physical Hilbert space
$\Hil^{(S)}$ of the PTSQM theory. The input Hamiltonian $\hat{H}$
may be declared self-adjoint in $\Hil^{(S)}$.  In this space it also
generates the correct unitary time-evolution of the system in
question, at least for the time-independent interactions, \cite{timedep}.

\section{Summary}

We have given here a brief review of some aspects of PTSQM with
particular interest to the role of different inner products and
their related conjugations. We have also discussed how the time
evolution of a system with a non self-adjoint Hamiltonian can be
analyzed within this settings, and the role of the different inner
products is discussed.

Among other lines of research, we believe that the construction of
algebras of unbounded operators associated to PTSQM, along the same
lines as in \cite{bagrev}, is an interesting task and we hope to be
able to do that in the near future.


\section*{Acknowledgment}

Work supported in part by the GA\v{C}R grant Nr. P203/11/1433, by
the M\v{S}MT ``Doppler Institute" project Nr. LC06002 and by the
Inst. Res. Plan AV0Z10480505 and in part by M.I.U.R.


\end{document}